# Limit Theorems for Differential Systems with Multipoint Boundary Conditions in Sobolev Spaces

**Olena Atlasiuk[1,3] · Vladimir Mikhailets[2,3]**

Received: 4 August 2026 / Accepted: 13 August 2026

## Abstract

The most general class of multipoint inhomogeneous boundary-value problems for systems of linear ordinary differential equations of arbitrary order $r \geq 1$ is investigated, whose solutions belong to a given Sobolev space $W_p^{n+r}$, where $n \geq 0$ and $1 \leq p \leq \infty$. The boundary conditions in these problems contain Caputo derivatives of fractional or integer orders, which may exceed the order of the equations of the differential system. Constructive sufficient conditions are established under which the solutions of these problems are continuous in the Sobolev space $W_p^{n+r}$ with respect to a parameter from an abstract metric space.

## 1 Introduction and Problem Statement

Multipoint boundary-value problems with derivatives whose order is lower than the order of the differential system/equation are classical objects of study in the theory of ordinary differential equations. Unlike Cauchy problems, solutions of such problems

✉ Olena Atlasiuk
olena.atlasiuk@helsinki.fi

Vladimir Mikhailets
vladimir.mikhailets@gmail.com

1 University of Helsinki, Department of Mathematics and Statistics, Pietari Kalmin katu 5, P.O. Box 68, 00014 Helsinki, Finland

2 King's College London, Strand, London WC2R 2LS, UK

3 Institute of Mathematics of the National Academy of Sciences of Ukraine, st. Tereschenkivska 3, 01024 Kyiv, Ukraine

may not exist or may not be unique. If the coefficients and the right hand side of the differential system are sufficiently smooth, then its solutions have smoothness exceeding the order of the system. Therefore, in this case, multipoint boundary conditions may contain derivatives of higher orders. The question of solvability of such boundary-value problems in Sobolev spaces was investigated in the paper [8]. Moreover, the index of the inhomogeneous boundary-value problem and necessary and sufficient conditions for its well-posedness were found there in given Sobolev spaces.

The issue of the continuous dependence of solutions of linear boundary-value problems on a numerical (continuous or discrete) parameter was actively studied in the articles of Kiguradze [9], Ashordia [1–3], Ponomaryov [10], Mikhailets and his students [4, 5]. The present paper continues this line of research. Its main features are as follows:

— For the first time, we study multipoint inhomogeneous boundary-value problems with fractional Caputo derivatives of different orders. The order of these derivatives may exceed the order of the differential system/equation.
— We consider and investigate for the first time boundary-value problems with the parameter $\mu$ on the left and right-hand sides, which is an element of an abstract metric space. This allows us to cover the cases of continuous, discrete and functional parameters from a single point of view, which is important for applications.
— We investigate the general case when the number of points and their localization depend on the parameter, which fundamentally complicates the analysis of the problems.
— We consider the case of arbitrary Sobolev spaces on a finite interval $W_p^n$, where $n \in \mathbb{N} \cup \{0\}$, $1 \leq p \leq \infty$.
— Not only limit theorems for solutions of boundary-value problems have been established, but also the estimate of their convergence in terms of the discrepancy of the solutions has been found.

Let a finite interval $(a, b) \subset \mathbb{R}$ be given and the parameters

$$\{m, r\} \subset \mathbb{N}, \quad n \in \mathbb{N} \cup \{0\}, \quad 1 \leq p \leq \infty.$$

We denote by $W_p^n := W_p^n([a, b]; \mathbb{C})$ the complex Sobolev space and set $W_p^0 := L_p$. We also denote by $(W_p^n)^m := W_p^n([a, b]; \mathbb{C}^m)$ and $(W_p^n)^{m\times m} := W_p^n([a, b]; \mathbb{C}^{m\times m})$ the Sobolev spaces of vector functions and matrix functions, respectively, whose elements belong to the functional space $W_p^n$. We denote the norms in these spaces by $\|\cdot\|_{n,p}$, they are the sums of the corresponding norms in $W_p^n$ of all elements of the vector or matrix function. It is always clear from the context to which space (scalar, vector or matrix functions) the norm refers. If $m = 1$, then all these spaces coincide. As is known, the spaces $W_p^n$ are Banach spaces with respect to the norm

$$\|y\|_{n+r,p} = \sum_{s=0}^{n+r} \|y^{(s)}(\cdot)\|_{L_p}.$$

They are separable if and only if $p < \infty$.

We investigate inhomogeneous boundary conditions with the parameter $\mu$ belonging to an abstract metric space $\mathcal{M}$. This allows us to study from a single point of view the cases when the parameter $\mu$ is discrete/continuous numerical, or functional.

Consider a system of linear differential equations of order $r$ with parameter:

$$L(\mu)y(\mu) := y^{(r)}(t,\mu)+\sum_{s=1}^{r} A_{r-s}(t,\mu)y^{(r-s)}(t,\mu) = f(t,\mu), \quad t \in (a,b). \tag{1.1}$$

Here, the matrix functions $A_{r-s}(\cdot,\mu) \in (W_p^n)^{m\times m}$, the right-hand sides $f(\cdot,\mu) \in (W_p^n)^m$, the unknown vector function $y(\cdot,\mu) \in (W_p^{n+r})^m$ satisfy eq. (1.1) at each point of the interval $(a,b)$ if $n \in \mathbb{N}$, and almost everywhere if $n = 0$. In addition, the solution $y(\cdot,\mu)$ should satisfy the inhomogeneous boundary condition

$$B(\mu)y(\cdot,\mu) = q(\mu) \in \mathbb{C}^{rm}, \tag{1.2}$$

where $B(\mu)$ is a linear continuous operator

$$B(\mu)\colon (W_p^{n+r})^m \to \mathbb{C}^{rm}.$$

The boundary-value problem (1.1), (1.2) is represented by the linear operator $\big(L(\mu), B(\mu)\big)\colon y \mapsto \big(L(\mu)y, B(\mu)y\big)$. Then, as proven in the paper [8, Theorem 2.1], for each fixed value of the parameter $\mu \in \mathcal{M}$, the linear operator $\big(L(\mu), B(\mu)\big)$ is a continuous Fredholm operator with index zero.

Let $\mu_0$ be an arbitrary limit point in the metric space $\mathcal{M}$. We shall assume the following throughout the paper.

**(0)** *The boundary condition* $\big(L(\mu_0), B(\mu_0)\big)y = 0$ *has only the trivial solution.*

Then every inhomogeneous boundary-value problem

$$\big(L(\mu_0), B(\mu_0)\big)y = (f,q) \in (W_p^{n+r})^m \oplus \mathbb{C}^{rm}$$

has a solution $y(\cdot,\mu_0)$ and it is unique. Now let $\mu \to \mu_0$ in $\mathcal{M}$. Then the following questions naturally arise:

1. What conditions should be satisfied on the coefficient matrices of the system and the boundary condition operators $B(\mu)$ to guarantee the invertibility of the operators $\big(L(\mu), B(\mu)\big)$, where $\mu$ belongs to some neighborhood of the point $\mu_0$? That is, under what conditions does the solution $y(\cdot,\mu)$ of the operator equation

$$\big(L(\mu), B(\mu)\big)y(\cdot,\mu) = \big(f(\cdot,\mu), q(\mu)\big)$$

   exist and is unique.
2. What conditions should be imposed on the left-hand sides of the problems to guarantee the following fact: if $f(\cdot,\mu) \to f(\cdot,\mu_0)$, $q(\mu) \to q(\mu_0)$, $\mu \to \mu_0$, then

$$\big\|y(\cdot,\mu) - y(\cdot,\mu_0)\big\|_{n+r,p} \to 0, \quad \text{as} \quad \mu \to \mu_0. \tag{1.3}$$

As established in the papers [6, 8], in the general case the condition (1.3) is equivalent to the fact that

$$L(\mu) \to L(\mu_0), \quad B(\mu) \to B(\mu_0), \quad \text{as} \quad \mu \to \mu_0$$

in strong operator topologies of Banach spaces.

In this article, we explore the case when the boundary conditions $B(\mu) = q(\mu)$ are multipoint. They may contain Caputo derivatives of fractional or integer orders exceeding the order of the differential system $r$.

In the case where all derivatives are classical and $\mathcal{M} = [0, \varepsilon]$, $\mu_0 = 0$, similar results were established in the works [4, 5].

## 2 Main Results

For the convenience of the reader, we present some definitions and well-known facts from fractional calculus, which we will use further.

Let $l \geq 0$. The right fractional Riemann-Liouville derivative of order $l$ of a function $y$ is defined as

$$\left(\mathrm{D}^l_{a+} y\right)(x) := \frac{1}{\Gamma(n-l)} \left(\frac{\mathrm{d}}{\mathrm{d}x}\right)^n \int_a^x \frac{y(t)\mathrm{d}t}{(x-t)^{l-n+1}},$$

where $n = [l]+1$, $x > a$, $[l]$ is the integer part of the number $l$ [7, Section 2, Paragraph 2.1].

The right fractional Caputo derivative $\left({}^C\mathrm{D}^l_{a+} y\right)(x)$ of order $l \geq 0$ on $[a, b]$ of a function is defined as

$$\left({}^C\mathrm{D}^l_{a+} y\right)(x) := \left(\mathrm{D}^l_{a+}\left(y(t) - \sum_{s=0}^{n-1} \frac{y^{(s)}(a)}{s!}(t-a)^s\right)\right)(x),$$

where $n = [l] + 1$, for $l \notin \mathbb{N}_0$ and $n = l$ for $l \in \mathbb{N}_0$ [7, Section 2, Section 2.4]. Therefore, if $l \in \mathbb{N}_0$, then the Caputo derivative coincides with the ordinary derivative.

If the function $y \in W_p^{n+r}$, then its fractional Caputo derivatives of order $l \in [0, n+r]$ are well-defined and ${}^C\mathrm{D}^l_{a+} y \in W_p^{n+r-l}$ if $1 \leq p < \infty$ and ${}^C\mathrm{D}^l_{a+} y \in C^{n+r-l}$ if $p = \infty$.

We arbitrarily choose $N$ different points $\{t_1, \ldots, t_N\} \subset [a, b]$ and consider the operator $B(\mu_0)$, which corresponds to multipoint boundary conditions with Caputo derivatives at these points

$$B(\mu_0)y(\cdot, \mu_0) = \sum_{k=1}^{M} \sum_{j=1}^{N} \beta_j^{(l_k)}(\cdot, \mu_0) \left({}^C\mathrm{D}^{l_k}_{a+} y\right)(t_j) = q(\mu_0), \tag{2.1}$$

where $M \in \mathbb{N}$, the matrices $\beta_j^{(l_k)} \in \mathbb{C}^{rm \times m}$, and the orders of the derivatives satisfy the inequality $l_k < n + r - 1/p$.

Consider the inhomogeneous boundary-value problem

$$L(\mu_0)y(\cdot,\mu_0)=f(\cdot,\mu_0),\quad B(\mu_0)y(\cdot,\mu_0)=q(\mu_0), \tag{2.2}$$

as the boundary-value problem at $\mu\to\mu_0$ for such a multipoint boundary-value problem with parameter $\mu\in\mathcal{M}$

$$L(\mu)y(\cdot,\mu)=f(\cdot,\mu),\quad B(\mu)y(\cdot,\mu)=q(\mu).$$

In the sequel, $L(\mu)$ is given by the formula (1.1), and the multipoint boundary conditions have the following form

$$B(\mu)y(\cdot,\mu)=\sum_{j=0}^{N}\sum_{i=1}^{\omega_j(\mu)}\sum_{k=1}^{M}\beta_{j,i}^{(l_k)}(\cdot,\mu)\left({}^{C}\mathrm{D}_{a+}^{l_k}y\right)\big(t_{j,i}(\mu)\big)=q(\mu),\quad \mu\neq\mu_0. \tag{2.3}$$

Here, for each fixed value of the parameter $\mu$, the vector function $y(\cdot,\mu)\in(W_p^{n+r})^m$ is unknown, the numbers $l_k$ do not depend on $\mu$, and the matrices $\beta_{j,i}^{(l_k)}(\cdot,\mu)\in\mathbb{C}^{rm\times m}$. In the case of $\mu\neq\mu_0$ on the closed interval $[a,b]$, we choose at least $N$ points $t_{j,i}(\mu)$, which are connected in an $N+1$ series in the following way: for each fixed $j\in\{1,\dots,N\}$, all points $t_{j,i}(\mu)$ must have a common limit $t_j$ at $\mu\to\mu_0$, and for points $t_{0,i}(\mu)$, this requirement will not be imposed. The zero series for each fixed $\mu\neq\mu_0$ contains $\omega_0(\mu)$ arbitrary different points of the interval $[a,b]$. The series with index $j\in\{1,\dots,N\}$ contains $\omega_j(\mu)$ different points. Note that there may not be a zero series. In the boundary condition (2.3), we use the repeated sum over indices $j$ and $i$. This is due to further assumptions on the values of the parameter $j$ and the behavior of points $t_{j,i}(\mu)$, as $\mu\to\mu_0$.

The boundary-value problem with boundary conditions (2.3) can be written as an operator equation

$$\big(L(\mu),B(\mu)\big)y(\cdot,\mu)=\big(f(\cdot,\mu),q(\mu)\big),$$

where the linear operator

$$\big(L(\mu),B(\mu)\big)\colon (W_p^{n+r})^m\to(W_p^n)^m\oplus\mathbb{C}^{rm}. \tag{2.4}$$

According to [8, Theorem 2.1], it is a bounded Fredholm operator with index zero.

The main result of the article is to establish constructive sufficient conditions under which the solution $y=y(\cdot,\mu)$ of the multipoint boundary-value problem (2.3) is continuous in the parameter $\mu$ for $\mu=\mu_0$ in the Sobolev space $(W_p^{n+r})^m$, $1\le p\le\infty$. That is, the solution $y(\cdot,\mu)$ exists, is unique, and satisfies the boundary relation (1.3).

In order for the problem under study to make sense, we will assume that the above condition **(0)** is satisfied. This implies that for $\mu=\mu_0$ the Fredholm operator (2.4) is an isomorphism. Therefore, the boundary-value problem (2.2) has a unique solution $y(\cdot,\mu_0)\in(W_p^{n+r})^m$ and is uniquely defined for arbitrarily chosen right-hand sides $f(\cdot,\mu_0)\in(W_p^n)^m$ and $q(\mu_0)\in\mathbb{C}^{rm}$.

We formulate the main results of the paper in the form of Theorems 2.1, 2.2, and 2.3 below. The proofs of these theorems are given in Sects. 3 and 4.

Let us consider the **limit conditions** as $\mu \to \mu_0$:

(I) $A_{r-s}(\cdot;\mu) \to A_{r-s}(\cdot;\mu_0)$ in the space $(W_p^n)^{m\times m}$ for each $s \in \{1,\dots,r\}$;
(II) $B(\mu)y \to B(\mu_0)y$ in the space $\mathbb{C}^{rm}$ for each $y \in (W_p^{n+r})^m$.
(III) $f(\cdot,\mu) \to f(\cdot,\mu_0)$ in $(W_p^n)^m$, $q(\mu) \to q(\mu_0)$ in $\mathbb{C}^{rm}$,

and **assumption** as $\mu \to \mu_0$:

($\alpha$) $t_{j,i}(\mu) \to t_j$ for each $j \in \{1,\dots,N\}$, and $i \in \{1,\dots,\omega_j(\mu)\}$;
($\beta$) $\sum_{i=1}^{\omega_j(\mu)} \beta_{j,i}^{(l_k)}(\cdot,\mu) \to \beta_j^{(l_k)}(\cdot,\mu_0)$ for each $j \in \{1,\dots,N\}$, and for each $k \in \{1,\dots,M\}$, for which $l_k \le n+r-1$;
($\gamma$) $\sum_{i=1}^{\omega_j(\mu)} \left\|\beta_{j,i}^{(l_k)}(\cdot,\mu)\right\| \left|t_{j,i}(\mu)-t_j\right| \to 0$ for each $j \in \{1,\dots,N\}$, $i \in \{1,\dots,\omega_j(\mu)\}$, and $k \in \{1,\dots,M\}$;
($\delta$) $\sum_{i=1}^{\omega_0(\mu)} \left\|\beta_{0,i}^{(l_k)}(\cdot,\mu)\right\| \to 0$ for each $i \in \{1,\dots,\omega_0(\mu)\}$, and $k \in \{1,\dots,M\}$.

Here, the square numerical matrices $\beta_{j,i}^{(l_k)}(\cdot,\mu) \in \mathbb{C}^{m\times m}$. Note that for the zero-series points $t_{0,i}(\mu)$ there are no restrictions on their behavior as $\mu \to \mu_0$.

In the conditions ($\gamma$) and ($\delta$), the expression $\|\cdot\|$ is the norm of a complex numerical matrix. This norm is equal to the sum of the modules of all elements of the matrix. Assumptions ($\beta$) and ($\gamma$) allow that the norms of the coefficients $\left\|\beta_{j,i}^{(l_k)}(\cdot,\mu)\right\|$ can grow indefinitely as $\mu \to \mu_0$, but not too quickly.

We formulate a limit theorem for solutions to the multipoint boundary-value problem (2.3) in the case $p=\infty$.

**Theorem 2.1** *Let the boundary-value problem* (2.3)*, with $p=\infty$, satisfy assumptions* ($\alpha$)*,* ($\beta$)*,* ($\gamma$)*,* ($\delta$)*. Then it obeys the boundary condition* (II).

*If, moreover, the conditions* (0), (I)*, and* (III) *are satisfied, then for sufficiently small $d(\mu,\mu_0)$, a solution exists, is unique, and satisfies the limit relation* (1.3).

Here, $d(\cdot,\cdot)$ is a metric in a metric space $\mathcal{M}$.

Let us now turn to the case $1 \le p < \infty$. To do this, we consider the following **assumptions** as $\mu \to \mu_0$:

($\gamma_p$) $\sum_{i=1}^{\omega_j(\mu)} \left\|\beta_{j,i}^{(l_k)}(\cdot,\mu)\right\| \left|t_{j,i}(\mu)-t_j\right|^{n+r-1-1/p} = O(1)$ for each $j \in \{1,\dots,N\}$, $i \in \{1,\dots,\omega_j(\mu)\}$, and $k$, for which $l_k \in (n+r-1-1/p, n+r-1/p)$;
($\gamma'$) $\sum_{i=1}^{\omega_j(\mu)} \left\|\beta_{j,i}^{(l_k)}(\cdot,\mu)\right\| \left|t_{j,i}(\mu)-t_j\right| \to 0$ for each $j \in \{1,\dots,N\}$, $i \in \{1,\dots,\omega_j(\mu)\}$, and for each $k$, for which $l_k \le n+r-1-1/p$.

Observe that the systems of conditions ($\alpha$), ($\beta$), ($\gamma$), ($\delta$) and ($\alpha$), ($\beta$), ($\gamma_p$), ($\gamma'$), ($\delta$) do not guarantee uniform convergence of continuous operators $B(\mu)$ from $(W_p^{n+r})^m$ to $B(\mu_0)$ in $\mathbb{C}^{rm}$ as $\mu \to \mu_0$. Therefore, Theorems 2.1 and 2.2 do not follow from general facts of the theory of linear operators.

We formulate a limit theorem for solutions of the multipoint boundary-value problem (2.3) in the case $1 \leq p < \infty$.

**Theorem 2.2** *Let the boundary-value problem* (2.3)*, with* $1 \leq p < \infty$*, satisfy assumptions* $(\alpha)$, $(\beta)$, $(\gamma_p)$, $(\gamma')$, $(\delta)$. *Then it obeys the boundary condition* (II).

*If, moreover, the conditions* (0), (I) *and* (III) *are satisfied, then for sufficiently small* $d(\mu, \mu_0)$*, a solution exists, is unique, and satisfies the limit relation* (1.3).

The following theorem establishes the estimate of convergence in the limit ratio (1.3). Denote by the discrepancy

$$\widetilde{d}_{n,p}(\mu) := \left\| L(\mu)y(\cdot;\mu_0) - f(\cdot;\mu) \right\|_{n,p} + \left\| B(\mu)y(\cdot;\mu_0) - q(\mu) \right\|_{\mathbb{C}^{rm}}.$$

**Theorem 2.3** *Let the conditions of Theorem 2.1 be satisfied if* $p = \infty$*, or Theorem 2.2 be satisfied if* $1 \leq p < \infty$. *Then there exist positive numbers* $\varepsilon$, $\gamma_1$ *and* $\gamma_2$ *such that*

$$\gamma_1 \, \widetilde{d}_{n,p}(\mu) \leq \left\| y(\cdot;\mu_0) - y(\cdot;\mu) \right\|_{n+r,p} \leq \gamma_2 \, \widetilde{d}_{n,p}(\mu),$$

*for all* $\mu$ *from some ball* $\mathcal{B}(\mu_0, \varepsilon)$ *in the metric space* $\mathcal{M}$. *Here,* $\varepsilon$, $\gamma_1$ *and* $\gamma_2$ *do not depend on* $y(\cdot;\mu_0)$ *or* $y(\cdot;\mu)$.

## 3 Proof of Theorem 2.1

***Proof*** We write the operator $B(\mu)$ as a finite sum of $N+1$ terms divided by series:

$$B(\mu)y(\cdot,\mu) = B_0(\mu)y(\cdot,\mu) + B_1(\mu)y(\cdot,\mu) + \ldots + B_N(\mu)y(\cdot,\mu), \tag{3.1}$$

where

$$B_j(\mu)y(\cdot,\mu) = \sum_{i=1}^{\omega_0(\mu)} \sum_{k=1}^{M} \beta_{j,i}^{(l_k)}(\cdot,\mu) \left( {}^{C}\mathrm{D}_{a+}^{(l_k)} y \right)(t_{j,i}(\mu)), \quad j \in \{0, 1, \ldots, N\}.$$

The summation over the index $i$ in the formula for the operator $B_j$ at $\mu_0$ disappears because for this value of the parameter the series with number 0 is absent, and all the others consist of only one point $t_j$.

Let us prove the limit relations as $\mu \to \mu_0$

$$B_0(\mu) \to 0, \tag{3.2}$$

$$B_j(\mu) \to B_j(\mu_0), \quad j \in \{1, \ldots N\} \tag{3.3}$$

in strong operator topology. The convergence (3.2) and (3.3) guarantee the fulfillment of the limit condition (II).

We first show the strong convergence of the operators $B_0(\mu)$ to zero, that is, the fulfillment of the relation (3.2). With regard to the condition ($\delta$) and inequalities $\left\|\left({}^{C}\mathrm{D}_{a+}^{(l_k)}y\right)(t_{0,i}(\mu))\right\| \leq \|y\|_{n+r,\infty}$, we obtain the inequality

$$\sum_{i=1}^{\omega_0(\mu)} \left\|\beta_{0,i}^{(l_k)}(\cdot,\mu)\right\| \left\|\left({}^{C}\mathrm{D}_{a+}^{(l_k)}y\right)(t_{0,i}(\mu))\right\| \leq \sum_{i=1}^{\omega_0(\mu)} \left\|\beta_{0,i}^{(l_k)}(\cdot,\mu)\right\| \|y\|_{n+r,\infty} \to 0$$

for all admissible values of the parameters $k$. These and all other limits in the proof are considered under the condition that $\mu \to \mu_0$.

For an arbitrary vector function $y(\cdot,\mu) \in (W_\infty^{n+r})^m$, we have

$$\begin{aligned}
&\left\|B_j(\mu)y(\cdot,\mu) - B_j(\mu_0)y(\cdot,\mu_0)\right\| \\
&= \left\|\sum_{i=1}^{\omega_j(\mu)}\sum_{k=1}^{M} \beta_{j,i}^{(l_k)}(\cdot,\mu)\left({}^{C}\mathrm{D}_{a+}^{(l_k)}y\right)(t_{j,i}(\mu)) - \sum_{k=1}^{M}\beta_j^{(l_k)}(\cdot,\mu_0)\left({}^{C}\mathrm{D}_{a+}^{(l_k)}y\right)(t_j)\right\| \\
&\leq \sum_{i=1}^{\omega_0(\mu)}\sum_{k=1}^{M}\left\|\beta_{0,i}^{(l_k)}(\cdot,\mu)\right\| \left\|\left({}^{C}\mathrm{D}_{a+}^{(l_k)}y\right)(t_{0,i}(\mu))\right\| \\
&\quad + \sum_{k=1}^{M}\left\|\sum_{i=1}^{\omega_j(\mu)} \beta_{j,i}^{(l_k)}(\cdot,\mu)\left({}^{C}\mathrm{D}_{a+}^{(l_k)}y\right)(t_{j,i}(\mu)) - \beta_j^{(l_k)}(\cdot,\mu_0)\left({}^{C}\mathrm{D}_{a+}^{(l_k)}y\right)(t_j)\right\|.
\end{aligned} \tag{3.4}$$

We explore the second term on the right-hand side of the formula (3.4). For arbitrary $j \in \{1,\ldots,N\}$ and $k \in \{1,\ldots,M\}$, we have

$$\begin{aligned}
&\left\|\sum_{i=1}^{\omega_j(\mu)} \beta_{j,i}^{(l_k)}(\cdot,\mu)\big({}^{C}\mathrm{D}_{a+}^{(l_k)}y\big)(t_{j,i}(\mu)) - \beta_j^{(l_k)}(\cdot,\mu_0)\big({}^{C}\mathrm{D}_{a+}^{(l_k)}y\big)(t_j)\right\| \\
&\quad = \Bigg\|\sum_{i=1}^{\omega_j(\mu)} \beta_{j,i}^{(l_k)}(\cdot,\mu)\big({}^{C}\mathrm{D}_{a+}^{(l_k)}y\big)(t_{j,i}(\mu)) - \sum_{i=1}^{\omega_j(\mu)} \beta_{j,i}^{(l_k)}(\cdot,\mu)\big({}^{C}\mathrm{D}_{a+}^{(l_k)}y\big)(t_j) \\
&\qquad + \sum_{i=1}^{\omega_j(\mu)} \beta_{j,i}^{(l_k)}(\cdot,\mu)\big({}^{C}\mathrm{D}_{a+}^{(l_k)}y\big)(t_j) - \beta_j^{(l_k)}(\cdot,\mu_0)\big({}^{C}\mathrm{D}_{a+}^{(l_k)}y\big)(t_j)\Bigg\| \\
&\quad \leq \left\|\sum_{i=1}^{\omega_j(\mu)} \beta_{j,i}^{(l_k)}(\cdot,\mu)\Big(\big({}^{C}\mathrm{D}_{a+}^{(l_k)}y\big)(t_{j,i}(\mu)) - \big({}^{C}\mathrm{D}_{a+}^{(l_k)}y\big)(t_j)\Big)\right\| \\
&\qquad + \left\|\Big(\sum_{i=1}^{\omega_j(\mu)} \beta_{j,i}^{(l_k)}(\cdot,\mu) - \beta_j^{(l_k)}\Big)\big({}^{C}\mathrm{D}_{a+}^{(l_k)}y\big)(t_j)\right\|
\end{aligned}$$

$$\leq \sum_{i=1}^{\omega_j(\mu)} \|\beta_{j,i}^{(l_k)}(\cdot,\mu)\| \left\| \left({}^C\mathrm{D}_{a+}^{(l_k)} y\right)(t_{j,i}(\mu)) - \left({}^C\mathrm{D}_{a+}^{(l_k)} y\right)(t_j) \right\|$$
$$+ \left\| \sum_{i=1}^{\omega_j(\mu)} \beta_{j,i}^{(l_k)}(\cdot,\mu) - \beta_j^{(l_k)} \right\| \cdot \|y\|_{n+r,\infty}. \tag{3.5}$$

Then, condition $(\beta)$ guarantees that the following convergence holds:

$$\left\| \sum_{i=1}^{\omega_j(\mu)} \beta_{j,i}^{(l_k)}(\cdot,\mu) - \beta_j^{(l_k)}(\cdot,\mu_0) \right\| \|y\|_{n+r,\infty} \to 0. \tag{3.6}$$

From the properties of the fractional differentiation operation and the inequality $n + r - l_k \geq 1$ the next result follows: if $y(\cdot,\mu) \in (W_\infty^{n+r})^m$, then the derivative ${}^C\mathrm{D}_{a+}^{(l_k)} y(\cdot,\mu)$ is a Lipschitz continuous function. Hence, the relation holds

$$\sum_{i=1}^{\omega_j(\mu)} \left\|\beta_{j,i}^{(l_k)}(\cdot,\mu)\right\| \left\| \left({}^C\mathrm{D}_{a+}^{(l_k)} y\right)(t_{j,i}(\mu)) - \left({}^C\mathrm{D}_{a+}^{(l_k)} y\right)(t_j) \right\| \to 0. \tag{3.7}$$

From formulas (3.2), (3.5), (3.6), (3.7) for $j \in \{1,\ldots,N\}$, the convergence directly follows

$$\left\| \sum_{i=1}^{\omega_j(\mu)} \beta_{j,i}^{(l_k)}(\cdot,\mu) \left({}^C\mathrm{D}_{a+}^{(l_k)} y\right)(t_{j,i}(\mu)) - \beta_j^{(l_k)}(\cdot,\mu_0) \left({}^C\mathrm{D}_{a+}^{(l_k)} y\right)(t_j) \right\| \to 0. \tag{3.8}$$

The convergence of (3.8) implies the fulfillment of the condition (3.3). Therefore, in view of (3.2) and (3.3), we conclude that $\|B(\mu)y(\cdot,\mu) - B(\mu_0)y(\cdot,\mu_0)\| \to 0$. But the vector function $y(\cdot,\mu) \in (W_\infty^{n+r})^m$ is arbitrary. Consequently, the limit condition (II) is fulfilled.

The first statement of Theorem 2.1 is proved. The second statement follows from the proof above and [6, Theorem 2.2]. □

## 4 Proof of Theorems 2.2 and 2.3

***Proof of Theorem 2.2*** As before, we assume that the operator $B(\mu)$ is written in the form (3.1), where the sum of $N + 1$ terms is divided by series of points.

Based on the Banach–Steinhaus theorem, it suffices to show that the norms of the operators $B(\mu)\colon (W_p^{n+r})^m \to \mathbb{C}^{rm}$ are bounded in some neighborhood of the point $\mu_0 \in \mathcal{M}$, and $B(\mu)y(\cdot,\mu) \to B(\mu_0)y(\cdot,\mu_0)$ as $\mu \to \mu_0$ for each vector function $y(\cdot,\mu)$ that belongs to the dense set of infinitely differentiable vector functions on the interval $[a,b]$ in the space $(W_p^{n+r})^m$.

We first prove the uniform boundedness of the norms of the operator in $\mu$

$$B(\mu) = \sum_{j=0}^{N} B_j(\mu).$$

We choose the parameter $\mu$ from a sufficiently small neighborhood of the point $\mu_0$ in $\mathcal{M}$ and an arbitrary vector function $y(\cdot,\mu) \in (W_p^{n+r})^m$. By virtue of the limit conditions (2.1) and (2.3), the inequality holds

$$\begin{aligned}
&\left\| B_j(\mu)y(\cdot,\mu) - B_j(\mu_0)y(\cdot,\mu_0) \right\| \\
&= \left\| \sum_{i=1}^{\omega_j(\mu)} \sum_{k=1}^{M} \beta_{j,i}^{(l_k)}(\cdot,\mu) \left({}^{C}\mathrm{D}_{a+}^{(l_k)} y\right)(t_{j,i}(\mu)) - \sum_{k=1}^{M} \beta_j^{(l_k)}(\cdot,\mu_0) \left({}^{C}\mathrm{D}_{a+}^{(l_k)} y\right)(t_j) \right\| \\
&\sum_{i=1}^{\omega_0(\mu)} \sum_{k=1}^{M} \left\| \beta_{0,i}^{(l_k)}(\cdot,\mu) \right\| \left\| \left({}^{C}\mathrm{D}_{a+}^{(l_k)} y\right)(t_{0,i}(\mu)) \right\| + \\
&\sum_{k=1}^{M} \left\| \sum_{i=1}^{\omega_j(\mu)} \beta_{j,i}^{(l_k)}(\cdot,\mu) \left({}^{C}\mathrm{D}_{a+}^{(l_k)} y\right)(t_{j,i}(\mu)) - \beta_j^{(l_k)}(\cdot,\mu_0) \left({}^{C}\mathrm{D}_{a+}^{(l_k)} y\right)(t_j) \right\|.
\end{aligned}$$

We show that the norm of the operator corresponding to the zero series is bounded. Using the continuity of the embedding of Sobolev spaces into Hölder spaces, we obtain the inequality

$$\sum_{i=1}^{\omega_0(\mu)} \left\| \beta_{0,i}^{(l_k)}(\cdot,\mu) \right\| \left\| \left({}^{C}\mathrm{D}_{a+}^{(l_k)} y\right)(t_{0,i}(\mu)) \right\| \le c_0 \sum_{i=1}^{\omega_0(\mu)} \left\| \beta_{0,i}^{(l_k)}(\cdot,\mu) \right\| \left\| y \right\|_{n+r,p} \tag{4.1}$$

for all admissible values of the indices $l_k$, where $c_0$ is the norm of the embedding operator.

Furthermore, for all $k \in \{1,\dots,M\}$ and $j \in \{1,\dots,N\}$, we have

$$\begin{aligned}
&\left\| \sum_{i=1}^{\omega_j(\mu)} \beta_{j,i}^{(l_k)}(\cdot,\mu)\big({}^{C}\mathrm{D}_{a+}^{(l_k)} y\big)(t_{j,i}(\mu)) - \beta_j^{(l_k)}(\cdot,\mu_0)\big({}^{C}\mathrm{D}_{a+}^{(l_k)} y\big)(t_j) \right\| \\
&\quad = \left\| \sum_{i=1}^{\omega_j(\mu)} \beta_{j,i}^{(l_k)}(\cdot,\mu)\big({}^{C}\mathrm{D}_{a+}^{(l_k)} y\big)(t_{j,i}(\mu)) - \sum_{i=1}^{\omega_j(\mu)} \beta_{j,i}^{(l_k)}(\cdot,\mu)\big({}^{C}\mathrm{D}_{a+}^{(l_k)} y\big)(t_j) \right. \\
&\qquad \left. + \sum_{i=1}^{\omega_j(\mu)} \beta_{j,i}^{(l_k)}(\cdot,\mu)\big({}^{C}\mathrm{D}_{a+}^{(l_k)} y\big)(t_j) - \beta_j^{(l_k)}(\cdot,\mu_0)\big({}^{C}\mathrm{D}_{a+}^{(l_k)} y\big)(t_j) \right\| \\
&\quad \le \left\| \sum_{i=1}^{\omega_j(\mu)} \beta_{j,i}^{(l_k)}(\cdot,\mu)\Big(\big({}^{C}\mathrm{D}_{a+}^{(l_k)} y\big)(t_{j,i}(\mu)) - \big({}^{C}\mathrm{D}_{a+}^{(l_k)} y\big)(t_j)\Big) \right\| \\
&\qquad + \left\| \Big( \sum_{i=1}^{\omega_j(\mu)} \beta_{j,i}^{(l_k)}(\cdot,\mu) - \beta_j^{(l_k)}(\cdot,\mu_0) \Big)\big({}^{C}\mathrm{D}_{a+}^{(l_k)} y\big)(t_j) \right\|
\end{aligned}$$

$$\leq \sum_{i=1}^{\omega_j(\mu)} \|\beta_{j,i}^{(l_k)}(\cdot,\mu)\| \cdot \|({}^C\mathrm{D}_{a+}^{(l_k)} y)(t_{j,i}(\mu)) - ({}^C\mathrm{D}_{a+}^{(l_k)} y)(t_j)\|$$
$$+ \left\| \sum_{i=1}^{\omega_j(\mu)} \beta_{j,i}^{(l_k)}(\cdot,\mu) - \beta_j^{(l_k)}(\cdot,\mu_0) \right\| \cdot \|y\|_{n+r,p}.$$

Hence, taking into account that ${}^C\mathrm{D}_{a+}^{(l_k)} : (W_p^{n+r})^m \to (W_p^{n+r-l_k})^m$ and the embedding theorem, we obtain the required boundedness of the operators $B(\mu)$ by the norm in some neighborhood of the point $\mu_0 \in \mathcal{M}$.

We now justify the strong convergence of the operators $B(\mu)$ to $B(\mu_0)$. In view of the condition $(\delta)$ and the inequality (4.1), we obtain

$$\sum_{i=1}^{\omega_0(\mu)} \left\|\beta_{0,i}^{(l_k)}(\cdot,\mu)\right\| \left\|\left({}^C\mathrm{D}_{a+}^{(l_k)} y\right)(t_{0,i}(\mu))\right\| \to 0. \tag{4.2}$$

In addition, under the condition $(\beta)$

$$c_0 \left\| \sum_{i=1}^{\omega_j(\mu)} \beta_{j,i}^{(l_k)}(\cdot,\mu) - \beta_j^{(l_k)}(\cdot,\mu_0) \right\| \|y\|_{n+r,p} \to 0. \tag{4.3}$$

If the vector function $y$ belongs to $(C^\infty)^m$, then

$$\begin{aligned} &\sum_{i=1}^{\omega_j(\mu)} \left\|\beta_{j,i}^{(l_k)}(\cdot,\mu)\right\| \left\|\left({}^C\mathrm{D}_{a+}^{(l_k)} y\right)(t_{j,i}(\mu)) - \left({}^C\mathrm{D}_{a+}^{(l_k)} y\right)(t_j)\right\| \\ &\leq \sum_{i=1}^{\omega_j(\mu)} \left\|\beta_{j,i}^{(l_k)}(\cdot,\mu)\right\| \max_{a\le t\le b} \left\|\left({}^C\mathrm{D}_{a+}^{(l_k+1)} y\right)(t)\right\| \left|t_{j,i}(\mu) - t_j\right| \to 0 \end{aligned} \tag{4.4}$$

for all $l_k$ having regard to the conditions $(\alpha)$ and $(\beta)$. Since the convergence of operators in a strong operator topology implies their uniform boundedness by the Banach–Steinhaus theorem, then from formulas (4.2), (4.3), (4.4) and the boundedness of the norms of the operators $B(\mu)$ in the neighborhood of the point $\mu_0 \in \mathcal{M}$, we obtain the convergence of $B(\mu)y(\cdot,\mu) \to B(\mu_0)y(\cdot,\mu_0)$ in $\mathbb{C}^{rm}$ as $\mu \to \mu_0$ for each $y(\cdot,\mu) \in (W_p^{n+r})^m$.

The first statement of Theorem 2.2 is proved. The second statement follows from the above and [6, Theorem 2.2]. □

The statement of Theorem 2.3 follows from the already proven Theorems 2.1, 2.2, and [6, Theorem 2.5].

**Acknowledgements** We thank the reviewer for the helpful comments, which allowed us to improve the presentation of the results of the paper.

**Author Contributions** All authors contributed equally to the study, read and approved the final version of the submitted manuscript.

**Funding** Open Access funding provided by University of Helsinki (including Helsinki University Central Hospital). The work of the first-named author was funded by Postdoctoral Fellowship EU-MSCA4Ukraine (number: 1244691, WBS-number: 4100609).
This project has received funding through the MSCA4Ukraine project, which is funded by the European Union. Views and opinions expressed are however those of the author only and do not necessarily reflect those of the European Union, the European Research Executive Agency or the MSCA4Ukraine Consortium. Neither the European Union nor the European Research Executive Agency, nor the MSCA4Ukraine Consortium as a whole nor any individual member institution of the MSCA4Ukraine Consortium can be held responsible for them.
The second named author would like to thank the Isaac Newton Institute for Mathematical Sciences, Cambridge, for support and hospitality during the "Solidarity Program" where work on this paper was undertaken. This work was supported by "EPSRC grant no EP/R014604/1". The author wishes to thank the Department of Mathematics, King’s College London, for their hospitality, and to the Ministry of Education and Science of Ukraine for support under the grant 0126U000898.

**Data Availability** No datasets were generated or analysed during the current study.

## Declarations

**Conflict of interest** The authors declare that they have no conflict of interest regarding the publication of this paper.

## References

1. Ashordia, M.: Conditions of the existence and uniqueness of solutions of the multipoint boundary value problem for a system of generalized ordinary differential equations. Georgian Math. J. **5**, 1–24 (1998). https://doi.org/10.1515/GMJ.1998.1
2. Ashordia, M.: On the general and multipoint boundary value problems for linear systems of generalized ordinary differential equations, linear impulse and linear difference systems. Mem. Differ. Equ. Math. Phys. **36**, 1–80 (2005)
3. Ashordia, M.: The general boundary value problems for linear systems of generalized ordinary differential equations, linear impulsive differential and ordinary differential systems. Numerical solvability. Mem. Differ. Equ. Math. Phys **81**, 1–184 (2020)
4. Atlasiuk, O.M.: Limit theorems for solutions of multipoint boundary-value problems in Sobolev spaces. J. Math. Sci. **247**, 238–247 (2020). https://doi.org/10.1007/s10958-020-04799-w
5. Atlasiuk, O.M.: Limit theorems for solutions of multipoint boundary-value problems with a parameter in Sobolev spaces. Ukrain. Math. J. **72**, 1175–1184 (2021). https://doi.org/10.1007/s11253-020-01859-x
6. Atlasiuk, O., Mikhailets, V., Taskinen, J.: Parameter-dependent inhomogeneous boundary-value problems in Sobolev spaces, Preprint (2025). https://doi.org/10.48550/arXiv.2512.21361
7. Kilbas, A.A., Srivastava, H.M., Trujillo, J.J.: Theory and Applications of Fractional Differential Equations, Elsevier, Amsterdam (2006). https://doi.org/10.1016/S0304-0208(06)80001-0
8. Mikhailets, V., Atlasiuk, O.: The solvability of inhomogeneous boundary-value problems in Sobolev spaces. Banach J. Math. Anal. (2024). https://doi.org/10.1007/s43037-023-00316-8
9. Kiguradze, I.T.: On boundary-value problems for linear differential systems with singularities. Differ. Equ. **39**, 212–225 (2003). https://doi.org/10.1023/A:1025152932174

10. Ponomarev, V.D.: Necessary and sufficient conditions for the solvability of a multipoint boundary value problem for ordinary differential equations of the first order. Differ. Equ. **14**, 661–663 (1978)

**Publisher's Note** Springer Nature remains neutral with regard to jurisdictional claims in published maps and institutional affiliations.